\documentclass[a4paper, 11pt]{amsart}
\pdfoutput=1

\usepackage{lmodern}
\usepackage{microtype}
\usepackage{amssymb, amsmath}
\usepackage[english]{babel}
\usepackage{graphicx}
\usepackage{url}
\usepackage{cite}
\usepackage{a4wide}
\usepackage{amscd}
\usepackage[pdftex, plainpages = false, pdfstartview=FitH,
    pdfpagelabels,  pdfpagelayout = SinglePage, 
    bookmarks = true, bookmarksopen = true, 
    bookmarksnumbered = true, linktocpage,
    colorlinks = true, 
    linkcolor = blue, urlcolor  = blue, 
    citecolor = red,  anchorcolor = green, 
    hyperindex = true, hyperfigures]{hyperref} 

\DeclareGraphicsExtensions{.png, .jpg, .pdf}
\graphicspath{{img/}{img/}{img/}}

\title[Toric arrangements]%
			{Deletion-restriction in toric arrangements}
\author[P\, Deshpande ]{Priyavrat Deshpande}
\author[K\, Sutar]{Kavita Sutar}
\address{Chennai Mathematical Institute \\ SIPCOT IT park, Siruseri\\ Tamil Nadu, India}
\email{pdeshpande@cmi.ac.in} 
\email{ksutar@cmi.ac.in}


\def\ds{\displaystyle}

\newcommand{\ol}{\overline}

\newcommand{\beq} {\begin{eqnarray}}
\newcommand{\eeq}{\end{eqnarray}}

\def\@begintheorem#1#2{\par\bgroup{\sc #1\ #2. }\it \ignorespaces}
\def\@opargbegintheorem#1#2#3{\par\bgroup{\sc #1\ #2\ (#3).}\it \ignorespaces}
\def\@endtheorem{\egroup}

\numberwithin{equation}{section}

\theoremstyle{plain}
\newtheorem{theorem}{Theorem}[section]
\newtheorem{lemma}[theorem]{Lemma}
\newtheorem{cor}[theorem]{Corollary}
\newtheorem{prop}[theorem]{Proposition}

\theoremstyle{definition}
\newtheorem{defn}[theorem]{Definition}
\newtheorem{ex}[theorem]{Example}
\newtheorem{xca}[theorem]{Exercise}
\newtheorem{conj}[theorem]{Conjecture}

\theoremstyle{remark}
\newtheorem{rem}[theorem]{Remark}

\newcommand{\bt}[1]{\begin{theorem}\label{#1}}
\newcommand{\bc}[1]{\begin{cor}\label{#1}}
\newcommand{\bl}[1]{\begin{lemma}\label{#1}}
\newcommand{\bp}[1]{\begin{prop}\label{#1}}
\newcommand{\be}[1]{\begin{ex}\label{#1}}
\newcommand{\bd}[1]{\begin{defn}\label{#1}}
\newcommand{\br}[1]{\begin{rem}\label{#1}}
\newcommand{\bx}[1]{\begin{xca}\label{#1}}
\newcommand{\bcon}[1]{\begin{conj}\label{#1}}

\newcommand{\et}{\end{theorem}}
\newcommand{\ec}{\end{cor}}
\newcommand{\el}{\end{lemma}}
\newcommand{\ep}{\end{prop}}
\newcommand{\ee}{\end{ex}}
\newcommand{\ed}{\end{defn}}
\newcommand{\exc}{\end{xca}}
\newcommand{\er}{\end{rem}}
\newcommand{\econ}{\end{conj}}

\newcommand{\bpr}{\begin{proof}}
\newcommand{\epr}{\end{proof}}

\newcommand{\ctor}[1]{T^{#1}_{\C}}

\def\A  {\mathcal{A}}

\def\C{\mathbb{C}}
\def\Z{\mathbb{Z}}

\begin{document}

\begin{abstract}
Deletion-restriction is a fundamental tool in the theory of hyperplane arrangements. Various important results in this field have been proved using deletion-restriction. In this paper we use deletion-restriction to identify a class of toric arrangements for which the cohomology algebra of the complement is generated in degree $1$. We also show that for these arrangements the complement is formal in the sense of Sullivan. 
\end{abstract}
 \maketitle
\section{Introduction}
A hyperplane arrangement is a finite collection $\A = \{H_1,\dots, H_n \}$ of hyperplanes in $\C^l$. 
One of the aspects of arrangement theory is to understand how the combinatorics of the intersections of these hyperplanes helps determine the topology of the complement of their union (denoted $M(\A)$). 
See for example \cite[Chapter 5]{orlik92} for an account of classical results that capture topological invariants of $M(\A)$ in terms of the associated combinatorics.
In this paper we focus on one such invariant: the cohomology algebra of the complement. 
We begin the introduction by recalling its structure. 

For every $1\leq p\leq n$ we say that a $p$-tuple $S = (H_1,\dots, H_p)$ of hyperplanes is \textit{independent} if $\dim(H_1\cap\cdots\cap H_p) = l - p$ and \textit{dependent} if the intersection is non-empty and of dimension strictly greater than $l-p$. 
Algebraically, the independence of $S$ means that the forms defining the hyperplanes in $S$ are linearly independent. 

Let $E$ be the free $\Z$-module generated by the elements $e_H$ for every $H\in \A$. 
Define $E(\A)$ to be the exterior algebra on $E$ and let $\partial$ denote the differential in $E(\A)$. For a $p$-tuple $S$ of hyperplanes we denote by $\bigcap S$ the intersection of elements in $S$ and by $e_S$ we mean $e_{H_1}\wedge\cdots\wedge e_{H_p}$. 
Let $I(\A)$ denote the ideal of $E(\A)$ generated by 
\[ \{ e_S~|~ \bigcap S = \emptyset\}\cup \{\partial e_S~|~ S\hbox{~is dependent}\}.\]

\bd{def15}
The \emph{Orlik-Solomon algebra} of a hyperplane arrangement $\A$ is the quotient algebra $E(\A)/I(\A)$ and is denoted by $A(\A)$.  \ed 

The following theorem combines the work of Arnol'd, Brieskorn, Orlik and Solomon (see \cite[Chapter 3, Section 5.4]{orlik92} for details).

\bt{thm1sec1} 
Let $\A = \{H_1, \dots, H_n\}$ be a hyperplane arrangement in $\C^l$. 
For $H\in \A$ choose a linear form $\alpha_H\in (\C^l)^*$, such that $\ker(\alpha_H) = H$. 
Then the integral cohomology algebra of the complement is generated by the classes represented by
\[\omega_H := \ds\frac{1}{2\pi i}d\log \alpha_H. \]
The map $\gamma\colon A(\A)\to H^*(M(\A); \Z)$ defined by \[\gamma(e_H)\mapsto \omega_H\]
induces an isomorphism of graded $\Z$-algebras.  \et 

There is a finer gradation of each $H^i(M(\A))$ indexed by the intersections. 
For $1\leq i\leq l$ let $W_{i_1},\dots, W_{i_k}$ denote all the codimension-$i$ subspaces that arise due to the intersections of the members of $\A$. 
For each $j = 1,\dots, k$ define 
\[M_{i_j} := \C^l\setminus \bigcup_{W_{i_j}\subseteq H} H. \]
$M_{i_j}$ is the complement of those hyperplanes in $\A$ that contain $W_{i_j}$.  
\bt{thm2sec1}
With notation as above, for each $1\leq i\leq l$, we have the following direct sum decomposition:
\[\ds H^i(M(\A); \Z) = \bigoplus_{j=1}^k H^i(M_{i_j}; \Z).\]
\et

For a distinguished hyperplane $H_0\in\A$ the arrangement $\A' := \A \setminus \{H_0\}$ is called the \textit{deleted arrangement} and $\A'' := \{H\cap H_0\mid H\in\A' \}$ is called the \textit{restricted arrangement}. Note that $\A''$ is an arrangement in $H_0 \cong \C^{l-1}$. 
The deletion-restriction property of hyperplane arrangements states that the long exact sequence in cohomology of the pair $(M(\A'), M(\A''))$ splits into short exact sequences 
\[0\to H^{k+1}(M(\A'))\to H^{k+1}(M(\A))\to H^k(M(\A''))\to 0. \]
We refer the reader to the article by Yuzvinski \cite{yuz01} for a more comprehensive account.
 
In this paper we consider toric arrangements, i.e., a finite collection of codimension-$1$ subtori in the complex $l$-torus $(\C^*)^l$. 
Our aim is to compute the cohomology algebra of the complement of the union of these subtori. 
To our knowledge the computation of such cohomology groups using the Leray spectral sequence were first performed by Looijenga in \cite{looijenga93}. 
Similar spectral sequence arguments have been used recently by Dupont in \cite{dupont13} in order to compute the Hodge decomposition of the complement of a divisor in a complex manifold and by Bibby in \cite{bibby13} where she computes cohomology groups of the complement of an abelian arrangement in an abelian variety. 
De Concini and Procesi used the theory D-modules for their cohomology computations in \cite{copr05, cpbook09}. 
They also compute the algebra structure in case of unimodular arrangements. 

In Theorem \ref{thm1} we find a sufficiency condition for the cohomology algebra to be generated in degree $1$. 
We use fairly elementary arguments to prove this result. 
Our treatment is based on the paper by Jozsa and Rice \cite{jozsarice91} where they give a simple reproof of the Brieskorn-Orlik-Solomon theorem using relative De Rham cohomology. We also remark that some of the results proved in this paper hold true for a more general situation as stated in \cite[Section 2]{jozsarice91}. 

Here is a quick sketch of the ideas we use. Let $X$ denote a finite-dimensional complex manifold and let $Y$ denote the union of finitely many codimension-$1$ submanifolds which intersect like hyperplanes. The cohomology groups of $M := X\setminus Y$, in principle, can be computed using the (cohomological) Gysin sequence which reads 
\begin{equation}
\cdots \to H^k(X)\stackrel{i^*}{\to} H^k(M)\stackrel{\mathrm{Res}}{\to} H^{k-1}(Y)\stackrel{\gamma}{\to} H^{k+1}(X)\to\cdots.
\end{equation} 
The map $i^*$ is induced by the inclusion $i\colon M\hookrightarrow X$. The group $H^{k-1}(Y)$ can be identified with the relative cohomology group $H^{k+1}(X, M)$ via the excision isomorphism and the Thom isomorphism. The map $\mathrm{Res}$ is the so-called \textit{Leray residue map} and $\gamma$ is the Gysin map. In our case $X$ is an $l$-torus and $Y$ is the union of finitely many codimension-$1$ subtori. We also point out to the reader that similar calculations were carried out by Sawyer in her thesis \cite{saw99}.\par 

The article is organized as follows. In Section \ref{sec2} we define toric arrangements. The cohomology calculations for the toric arrangements are carried out in Section \ref{sec3}. In particular we identify a class of toric arrangements, called deletion-restriction type, for which the cohomology algebra of the complement is generated in degree $1$. In Section \ref{sec4} we prove that for deletion-restriction type arrangements the complement is formal in the sense of Sullivan. We end the paper in Section \ref{sec5} by outlining future research. 

\section{Toric Arrangements}\label{sec2}
The standard $l$-dimensional complex torus is the space $(\C^*)^l$ of $l$-tuples of nonzero complex numbers. 
The torus is a group under multiplication of coordinates. 
In fact, it is an affine algebraic variety with the ring of \emph{Laurent polynomials} $\C[z_1^{\pm 1},\dots, z_l^{\pm 1}]$ as its coordinate ring. 

\bd{def1}
The \textit{character} of a torus is a multiplicative homomorphism $\chi\colon (\C^*)^l\to \C^*$ given by the evaluation of Laurent monomials  
\[\chi(z_1, \dots, z_l) = z_1^{n_1}\cdots z_l^{n_l}, \quad n_i\in \Z, \forall i. \] \ed

The following are some well-known facts regarding tori (see \cite[Section 5.2]{cpbook09}). 
The set $\Lambda$ of characters of $(\C^*)^l$ is a free abelian group isomorphic to $\Z^l$. 
Conversely, for a finitely generated abelian group $\Lambda$ of rank $l$ the variety $T^l_{\C} := \mathrm{Hom}(\Lambda, \C^*)$ is isomorphic to the product of an $l$-torus and a finite abelian group isomorphic to the torsion subgroup of $\Lambda$.

If $U_{\C}$ is the vector space of complex linear functionals on $\Lambda$ then for every $\phi\in U_{\C}$ the function $a\mapsto \exp(\phi(a))$ is a character. The connected, topologically closed subgroups (called the toric subgroups) of $\ctor{l}$ are isomorphic to $k$-tori for some $k\leq l$. In general a closed subgroup  of $\ctor{l}$ is isomorphic to $\ctor{k}\times A$ where $A$ is a finite abelian group. Any coset of a toric subgroup is homeomorphic to the group; thus topologically it is a torus. Let $W$ be a closed subgroup of $\ctor{l}$ and $W_0$ be the connected component containing identity. Then $W_0$ is a toric subgroup, the quotient $W/W_0$ is a finite abelian subgroup and $W$ is the union of the $|W/W_0|$ distinct cosets. \par 

Given a closed a subgroup $W$ of $\ctor{l}$ its inverse image under the exponential map $\exp: U_{\C}\to \ctor{l}$ is a closed subgroup of $U_{\C}$. In fact, there is a one-to-one correspondence between closed subgroups of $\ctor{l}$ and subgroups of $\Lambda$. Such subgroups are determined by integer matrices. An integer matrix $A$ of size $m\times l$ determines a mapping from $\ctor{l}$ to $\ctor{m}$. The kernel $W$ of this mapping is a closed subgroup of $\ctor{l}$ and every subgroup arises in this manner. Consequently $W$ depends only on the subgroup (i.e., the sub-lattice) of $\Lambda$ generated by the rows of $A$. Without loss of generality one can assume that rows of $A$ furnish a basis for this sub-lattice so that $m\leq l$.\par 

Given a character $\chi\in\Lambda$ and a non-zero complex number $c$ a \textit{toric hypersurface} $K_{\chi, c}$ is defined as the level set of $\chi$, i.e., $K_{\chi, c} := \{z\in (\C^*)^l\mid \chi(z) = c \}$. A toric hypersurface is a translate of a toric subgroup of codimension-$1$.

\bd{def2}
A toric arrangement is a finite collection 
\[\A := \{K_{\chi_1, c_1}, \dots,  K_{\chi_n, c_n}\mid \chi_i \in \Lambda, c_i\in \C^* ~\hbox{for}~ 1\leq i\leq n\}\]
of toric hypersurfaces in $(\C^*)^l$.
\ed

For notational simplicity we write $K_i$ instead of $K_{\chi_i, c_i}$. 
Without loss of generality we assume that each toric hypersurface in $\A$ is connected, i.e., each character is primitive (recall that a character is primitive if all the exponents are relatively prime). 
A toric arrangement can also be defined as the finite collection of pairs $(\chi, c)$ of a character and a non-zero complex number. 
Moreover, as each character is a Laurent monomial it is sometimes convenient to encode the arrangement information as the pair $(A, \overline{c})$, where $A$ is an $n\times l$ matrix whose $i$th row corresponds to exponents of $\chi_i$ and $\overline{c} = (c_1,\dots, c_n)$. We will use either of these notations as per convenience. 
To every toric arrangement $\A$ there is an associated periodic hyperplane arrangement $\tilde{\A}$ in $\C^l$. 
The hyperplanes in $\tilde{\A}$ are the inverse images $\exp^{-1}(K_i)$ for $1\leq i\leq n$. 
Note that the inverse image is the union of parallel integer translates of a codimension-$1$ subspace.\par 

\begin{defn}\label{intersection poset}
The \textit{intersection poset} $L(\A)$ of a toric arrangement $\A$ is the set of all connected components of all intersections of the toric hypersurfaces in $\A$ ordered by inclusion.
The elements of $L(\A)$ are called \textit{components of the arrangement}.
\end{defn}

The intersection poset is a ranked poset; the rank of every element is the dimension of the corresponding intersection.

Let $W$ be a component of a toric arrangement $\A$. 
For any $z\in W$ let $\A^z_W$ denote the hyperplane arrangement in a coordinate neighbourhood of $z$. 
More intrinsically, $\A^z_W = \{T_z(K)\mid W\subseteq K \}$ in the tangent space $T_z(\C^*)^l \cong \C^l$. 
Note that if $z$ is a point such that it does not belong to any other component contained in $W$ then the intersection data of $\A^z_W$ is independent of the choice of $z$. Hence one can safely disregard the reference to $z$ since we are mainly interested in the combinatorics of the arrangement. 
We refer to $\A_W$ as the \textit{local arrangement} at $W$. 
The complement $M(\A_W)$ is called the \textit{local complement} at $W$.\par 

An arrangement $\A = (A, \ol{c})$ is said to be \textit{unimodular} if it satisfies any one of the following equivalent conditions. 
\begin{enumerate}
\item Any intersection of members of $\A$ is either empty or connected. 
\item Any subset of rows of $A$ spans a direct summand in $\Lambda$.
\item All the $l\times l$ minors are $-1, 0$ or $1$ (i.e., $A$ is a unimodular matrix). 
\end{enumerate}

\bd{def3}
The \textit{complement} of a toric arrangement $\A$ in $(\C^*)^l$ is defined as follows 
\[M(\A) := (\C^*)^l \setminus \bigcup_{i=1}^{n} K_i.\]
\ed

\be{ex1} An arrangement in $\C^*$ is just a collection of finitely many points. The complement in this case has the homotopy type of wedge of circles.\ee

\be{ex2} A braid arrangement in $(\C^*)^l$ is the collection of $\binom{l}{2}$ toric hypersurfaces given by the following equations 
\[\A = \{z_i z_j ^{-1} = 1 \mid 1\leq i < j \leq l\}. \]
The complement in this case is the configuration space of $l$ ordered points in $\C^*$.
\ee

\section{Cohomology calculations}\label{sec3}

In this section we compute the de Rham cohomology of the toric complement $M(\A)$ for those arrangements $\A$ which satisfy a certain `deletion-restriction' criterion. As stated in the introduction the cohomology groups were computed by Looijenga \cite[2.4]{looijenga93} and independently by De Concini and Procesi in \cite[Theorem 4.2]{copr05}. We state their theorem below.
\bt{thm0} For every integer $0\leq k\leq l$ there is a non-canonical decomposition of the cohomology groups 
\[\ds H^k(M(\A)) = \bigoplus_{i = 0}^k [\bigoplus_{\mathrm{rank} W = i} H^{k - i}(W)\otimes H^i(M(\A_W))]\]
where $W$ is a component of the arrangement and $\A_W$ is the local arrangement at $W$.
\et 
Further, De Concini and Procesi also show that if $\A$ is a unimodular toric arrangement then $H^*_{{dR}}(M(\A))$ is generated in degree $1$ by logarithmic differential forms. They also compute the algebra structure and prove that the complement is formal in the sense of Sullivan \cite[Theorem 5.2]{copr05}.\par 

In this section we will show that there is a larger class of toric arrangements (containing the unimodular ones) for which the cohomology algebra of the complement is generated in degree $1$. 
We begin with an example of a non-unimodular arrangement such that the cohomology of its complement is generated in degree $1$.  

\be{ex4}Consider the arrangement $\A = \{z_1 = 1, z_2 = 1, z_1 z_2 = 1, z_1z_2^{-1} = 1 \}$ which is a non-unimodular arrangement as the intersection of the hypersurfaces corresponding to the last two characters is disconnected (it consists of two points $\{(1, 1),  (-1, -1)\}$). The complement of this arrangement is:
\[M(\A) = \{(z_1, z_2)\in\C^2 \mid z_1\neq 0, 1, z_2\neq 0, 1, z_1 z_2\neq 1, z_1z^{-1}_2 \neq 1\}. \]
Recall that $z\mapsto \frac{1 + z}{1 - z}$ is a bi-holomorphism from $\C\setminus \{0, 1\}$ to $\C\setminus \{\pm 1\}$. Applying this map coordinate-wise to $M(\A)$ we get 
\[M(\A) \cong \{(w_1, w_2)\in\C^2\mid w_1\neq\pm 1, w_2\neq\pm 1, w_1\pm w_2\neq 1\}. \]
The right hand side is a hyperplane complement hence $H^*(M(\A))$ is generated in degree $1$.  \ee

For notational simplicity we denote by $H^i(M)$ the de Rham cohomology of $M$ with complex coefficients. 
Let $\A = \{K_1,\dots, K_n\}$ be a toric arrangement. 
Our aim is to compute the cohomology of $M(\A)$; we do this by induction on $|\A|$. 
For $1\leq i\leq n$ we denote by $\A_i$ the sub-arrangement containing only the first $i$ subtori; we write  $M_i$ for $M(\A_i)$. 
By $\A_0$ we mean the empty arrangement and hence $M_0$ stands for $(\C^*)^l$.

Assume that the cohomology algebra of $M' := M_{i-1}$ is known for some $1\leq i\leq n$. Let $M''$ denote the intersection $M' \cap K_i$. Again for notational simplicity we denote the function $\chi_i - c_i$ by $f$. Hence $K_i = f^{-1}(0)$ and $0$ is a regular value of $f$. By the submersion theorem there is a tubular neighbourhood $U$ of $M''$ in $M'$. We consider $U$ as a rank $2$ trivial vector bundle and denote by $\pi\colon U\stackrel{}{\to} M''$ the projection map. This tubular neighbourhood is diffeomorphic to $M''\times \C$ via the map $\pi\times f|_U$. Denote the complement of the zero section by $U_0$.
The restricted map $\ol{\pi}\times \ol{f|_U}\colon U_0\to M''\times \C^*$ is a homotopy equivalence. By the K\"unneth isomorphism theorem we deduce that $H^*(U_0)\cong H^*(M'')\otimes H^*(\C^*)$ and hence 
\begin{equation}
H^k(U_0) = \{\ol{\pi}^*(\omega) + \ol{\pi}^*(\theta)\wedge (f|_{U_0})^*(\frac{dz}{z})\mid \omega\in H^k(M''), \theta\in H^{k-1}(M''), \frac{dz}{z}\in H^1(\C^*)\}\label{sec2eq1}.
\end{equation}

It follows from the Poincar\'e lemma that the inclusion induced homomorphism $H^*(U)\to H^*(U_0)$ is injective. Consequently the long exact sequence in cohomology for the pair $(U, U_0)$ breaks up into following short exact sequences: 
\begin{equation}
0 \to H^k(U) \hookrightarrow H^k(U_0)\to H^{k+1}(U, U_0)\to 0.\label{sec2eq2}
\end{equation}

\br{rem2} The relative cohomology group $H^k(U, U_0)$ can also be looked upon as the $k$-th local cohomology group of $U$ with support in $M''$ and coefficients in $\C$ (considered as the constant sheaf) and denoted $H^k_{M''}(U; \C)$. 
The group $H^{k+1}(U, U_0)$ is spanned by the classes represented by $\pi^*(\omega)\wedge d\log(f|_U)$ where $\omega$ ranges over representatives of $H^{k-1}(M'')$ and $d\log(f|_U) = (f|_U)^*(\frac{dz}{z})$. 
Since $U\stackrel{\pi}{\to} M''$ is trivial vector bundle of rank $2$ the cohomology of the pair $(U, U_0)$ is same as that of the associated Thom space. 
The Thom class in $H^2(U, U_0)$ is represented by $d\log(f|_U)$. \er

Let $M := M_i = M'\setminus M''$, whose cohomology we want to compute. 
The inclusion of pairs 

\[ \begin{CD}
M          @>\iota_M>>  M'\\
@AA{\ol{j}}A @AA{j}A\\
U_0 @>\iota_U>> U
\end{CD}\]

induces the following commuting diagram of long exact sequences of pairs

\begin{equation}
\begin{CD}
\cdots @>>> H^k(M')   @>{\iota^*_M}>>  H^k(M) @>>>   H^{k+1}(M', M) @>>> \cdots\\
@. @VV{j_k^*}V @VV{\ol{j_k}^*}V @V{\cong}V{\hbox{excision}}V @.\\
0  @>>> H^k(U) @>{\iota^*_U}>> H^k(U_0) @>>> H^{k+1}(U, U_0) @>>> 0.
\end{CD}\label{sec2eq3}
\end{equation}

Using Thom isomorphism we may replace $H^{k+1}(M', M)$ by $H^{k-1}(M'')$. 
We state a condition under which the top row of the above commuting diagram splits into short exact sequences; it resembles the deletion-restriction sequence described in the introduction. 

\bl{lem1}With notation as above if the map $j_k^*\colon H^k(M')\to H^k(U)$ is surjective then $H^*(M)$ is generated as a $\C$-algebra by $\iota_M^*(H^*(M'))$ and the class represented by $d\log(f|_M)$. \el

\bpr
Choose a class $\omega\in H^k(M'')$. From Equation \ref{sec2eq1} we have that $\ol{\pi}^*(\omega)\in H^k(U_0)$. As $j_k^*$ is surjective there exists $\alpha\in H^k(M')$ such that $j_k^*(\alpha) = \pi^*(\omega)$. Since $\ol{\pi} = \pi\circ\iota_U$ we have the following:
\begin{align}
\ol{\pi}^*(\omega) &= (\iota^*_U\circ \pi^*)(\omega)\nonumber \\
				   &= (\iota^*_U\circ j_k^*)(\alpha)\nonumber \\
				   &= (\ol{j_k}^*\circ \iota^*_M)(\alpha)\label{sec2eq4}. 
\end{align}
Hence $\ol{\pi}^*(\omega)$ is in the image of $\ol{j_k}^*$. 
Similarly one can prove that for $\theta\in H^{k-1}(M'')$ the class $\ol{\pi}^*(\theta) = \ol{j_k}^*(\iota^*_M(\beta))$ for some class $\beta\in H^{k-1}(M')$.
Using the identity $f|_{U_0} = f|_{M}\circ \ol{j_k}$ we get
\[\ol{\pi}^*(\theta)\wedge d\log(f|_{U_0}) = \ol{j_k}(\iota^*_M(\beta)\wedge d\log(f|_M)). \]
Therefore the map $\ol{j_k}$ is surjective. 
By commutativity of Diagram \ref{sec2eq3} one concludes that the connecting homomorphism in the top row is also surjective. 
As a result the top row breaks into the following deletion-restriction short exact sequence:
\begin{equation}
0\to H^k(M') \stackrel{\iota^*_M}{\hookrightarrow} H^k(M)\to H^{k-1}(M'')\to 0.\label{sec2eq5}
\end{equation}
For every $k\geq 0$, the group $H^k(M)$ is spanned by the following classes
\[\{\iota^*_M(\alpha)\mid \alpha\in H^k(M') \}\cup \{\iota^*_M(\beta)\wedge d\log(f|_M) \mid \beta\in H^{k-1}(M')\} \]
which proves the theorem. \epr

There are two inclusions $M_{i-1}\cap K_i\hookrightarrow U\stackrel{j}{\hookrightarrow} M_{i-1}$ of which the first inclusion is the homotopy inverse of the projection map. 
They induce the following maps in cohomology:
\[H^*(M_{i-1})\stackrel{j^*}{\to} H^*(U)\stackrel{\cong}{\to} H^*(M_{i-1}\cap K_i). \]
By abuse of notation we denote the above composition by $j^*$. 

\bt{thm1}
Let $\A$ be a toric arrangement in $(\C^*)^l$ and $M(\A)$ be its complement. If for every $1\leq i\leq n$ the inclusion induced map  $j^*\colon H^*(M_{i-1})\to H^*(M_{i-1}\cap K_i)$ is surjective then the algebra $H^*(M(\A))$ is generated by the classes of the logarithmic $1$-forms
\[\{d\log z_1, \cdots, d\log z_l, d\log (\chi_1-c_1),\cdots, d\log(\chi_n-c_n) \}. \] \et 

\bpr For $i=1$, $M_0 = (\C^*)^l$ and $M_0\cap K_1 = K_1$. The map $j^*$ is clearly surjective in this case. Lemma \ref{lem1} states that $H^*(M_0\setminus K_1) = H^*(M_1)$ is generated by $H^*(M_0)$ and $\{d\log(\chi_1-c_1)|_{M_1} \}$. Repeated application of Lemma \ref{lem1} proves the theorem. \epr

The above theorem gives a sufficiency condition for $H^*(M(\A))$ to be generated in degree $1$. 

\bt{thm2}
With notation as above, for all pairs $(M_{i-1}, M_{i-1}\cap K_i)$ for $1\leq i\leq n$ the following are equivalent:
\begin{enumerate}
\item The map $j^*\colon H^*(M_{i-1})\to H^*(M_{i-1}\cap K_i)$ is surjective.
\item The algebra $H^*(M_{i-1}\cap K_i)$ is generated by the classes represented by 
\begin{multline*}
\{
d\log z_1|_{M_{i-1}\cap K_i},\dots, d\log z_l|_{M_{i-1}\cap K_i},\\ 
d\log(\chi_1 - c_1)|_{M_{i-1}\cap K_i},\dots, d\log(\chi_{i-1} - c_{i-1})|_{M_{i-1}\cap K_i}
\}.
\end{multline*}
\end{enumerate}\et 

\bpr The second statement follows from the first using the proof of Theorem \ref{thm1}. 

For the converse we proceed by induction on $|\A|$. 
The case $i = 1$ is straightforward. Assume that the condition holds true for some $i-1$. Note that $K_r\cap M_{i-1} =\emptyset$ for all $1\leq r\leq i-1$. 
Consequently one can define a map from $g\colon M_{i-1}\to (\C^*)^{l+i-1}$ where the first $l$ coordinates of $g$ are restrictions of coordinate functions and the remaining $r$ coordinates  are the restrictions of $\chi_r-c_r$ to $M_{i-1}$ for $1\leq r\leq i-1$. 
For the same reasons we can define a map $h\colon M_{i-1}\cap K_i\to (\C^*)^{l+i-1}$ using the restrictions of the coordinates of $(\C^*)^l$ and of characters in $\A$. 
Now the result follows since in the commuting diagram below the maps $g^*$ and $h^*$ are surjective.
\begin{equation}
\begin{CD}
H^k(M_{i-1})   @>{j^*}>> H^k(U)\\
@A{g^*}AA @A{\cong}A{\pi^*}A      \\
H^{k}((\C^*)^{l+i-1})  @>{h^*}>> H^k(M_{i-1}\cap K_i).
\end{CD}\qedhere\label{sec2eq6}
\end{equation}
\epr

Note that the space $M_{i-1}\cap K_i$ is the complement of the toric arrangement in $K_i$ given by the toric hypersurfaces $\{K_r\cap K_i\mid 1\leq r\leq i-1 \}$. 
Hence it follows from the above theorem that the number of distinct connected components of $K_1\cap K_i,\dots, K_{i-1}\cap K_i$ is at most $i-1$. 

\be{ex5} Here we consider the arrangement $\A = \{z_1 = 1, z_2 = 1, z_1z_2 = 1, z_1z_2^{-1} = 1\}$ from Example \ref{ex4}. Since $M_0 = (\C^*)^2$ and $M_0\cap K_1 = K_1$ it is clear that $H^*(M_1)$ is generated by the restrictions of the classes $\{d\log z_1, d\log z_2, d\log(z_1-1) \}$. In fact, $M_1\cong (\C\setminus\{0, 1\})\times \C^*$.\par 
Now $M_1\cap K_2 = \{(z_1, 1)\in\C^2\mid z_1\neq 0, 1\}$ hence the cohomology $H^*(M_1\cap K_2)$ is generated by the restrictions of $\{d\log z_1, d\log(z_1-1) \}$. The case $M_2\cap K_3$ is similar.\par 
The cohomology of the complement $M_3\cap K_4$ is generated by $\{d\log z_1, d\log(z_1 - 1), d\log(z_1+1) \}$. Hence $H^*(M(\A))$ is generated by the $6$ classes represented by 
\[\{d\log z_1, d\log z_2, d\log(z_1 - 1), d\log(z_2 - 1), d\log(z_1z_2 - 1), d\log(z_1z_2^{-1} - 1) \}. \]
\ee

We now consider a non-example.

\be{ex6} Let $\A = \{z_1 = 1, z_1 z_2^3 = 1 \}$ in $(\C^*)^2$. Observe that $M_1\cap K_2 = \{(z_2^{-3}, z_2)\in \C^2 \mid z_2^3 \neq 0, 1\}$ which clearly means that its cohomology can not be generated by the restrictions of $\{d\log z_1, d\log z_2, d\log (z_1 - 1) \}$. Hence the hypothesis of Theorem \ref{thm1} is not satisfied. \ee

\bd{def4} 
A toric arrangement $\A$ is said to be of \emph{deletion-restriction type} if there exists an ordering $\{K_1,\dots, K_n\}$ on the hypersurfaces such that either of the conditions mentioned in Theorem \ref{thm2} holds true. 
Equivalently, for each $2\leq i\leq n$ the number of distinct connected components of $K_1\cap K_i,\dots, K_{i-1}\cap K_i$ is at most $i-1$.
\ed

The following result follows from the definitions.

\bp{prop1} A unimodular toric arrangement is of deletion-restriction type. \ep 

\section{Formality}\label{sec4}
An important property of the hyperplane complements is that they are formal in the sense of Sullivan. 
A consequence of Brieskorn's result stated in the introduction is that the cohomology algebra of the complement is isomorphic to the sub-algebra of rational differential forms generated by $1$ and the logarithmic forms $d\log\alpha_i$ ($\alpha_i$ is the linear form defining the $i$-th hyperplane). 
The aim of the current section is to extend this result to deletion-restriction type (DR-type for short) toric arrangements. 
We follow the strategy described in \cite[Section 3.5]{orlik92} and use the Leray's residue theory for differential forms. 
The main reference is the book by Pham \cite[Chapter III]{pham11}.\par 

Let $\A = \{K_1,\dots, K_n \}$ be a DR-type toric arrangement and let $\chi_i- c_i$ be the character defining $K_i$. The complement $M(\A)$ is an affine variety with the coordinate ring 
\[\ds S(\A) = \C[\Lambda]\left[\frac{1}{\prod_{i=1}^n (\chi_i-c_i)}\right].\] 
We denote by $(\Omega^{\bullet}, d)$ the algebraic de Rham complex of $M(\A)$ ($d$ is the usual de Rham differential). The algebra $\Omega^{\bullet}$ is a graded algebra explicitly described as 
\[\Omega^{\bullet} = S(\A)\otimes \wedge\langle \frac{dz_1}{z_1},\cdots, \frac{dz_l}{z_l}\rangle \]
where $S(\A)$ is assigned degree zero and the second part is the algebra of closed differential forms on $(\C^*)^l$. The grading on $\Omega^{\bullet}$ is given by
\[\ds\Omega^p = \bigoplus_{1\leq i_1<\cdots<i_p\leq l} S(\A) \frac{dz_{i_1}}{z_{i_1}}\wedge\cdots\wedge\frac{dz_{i_p}}{z_{i_p}}.\]

For notational simplicity we let $\xi_i = d\log z_i$ for $1\leq i\leq l$ and $\psi_j = d\log (\chi_j-c_j)$ for $1\leq j\leq n$. 
Let $R(\A)$ be the subalgebra of $\Omega^{\bullet}$ generated by $1$, the $1$-forms $\xi_1, \dots, \xi_l$ and $\psi_1,\dots, \psi_n$. Note that $R(\A)$ inherits the grading from $\Omega^{\bullet}$; we have $R_k(\A) := R(\A)\cap\Omega^k$ for $k\geq 0$. 
Similarly we define algebras $R(\A')$ and $R(\A'')$ where $\A' = \A_{i-1}$ for some $1\leq i \leq n$ ($\A_0$ is the empty arrangement) and $\A'' = \{K_1\cap K_i,\dots, K_{i-1}\cap K_i \}$. 
It is clear that there is an inclusion $\iota\colon R(\A')\hookrightarrow R(\A_i)$. \par 

We now define a homomorphism $\mathrm{res}: R(\A_i)\to R(\A'')$ using residues of differential forms such that $\mathrm{res}\circ\iota =0$. 
\bd{def5}
A differential form $\phi\in R(\A_i)$ is said to have at most a \textit{simple pole} along $M''$ if $(\chi_i-c_i)\wedge \phi$ is a form on $M'$. \ed
The following theorem guarantees that there are sufficient number of forms with simple poles along $M''$ (see \cite[Chapter III, Theorem 3.1]{pham11})

\bt{thm3}
Every closed regular form $\phi\in R(\A_i)$ is cohomologous to $\tilde{\phi}\in R(\A_i)$ which has at most a simple pole along $M''$. \et 

\bl{lem4}
Let $\phi\in R(\A_i)$ have at most a simple pole along $M''$. 
Then there exist regular forms $\ol{\phi}, \theta$ such that
\[\phi = \psi_i\wedge \ol{\phi} + \theta. \]
The forms $\ol{\phi}$ and $\theta$ do not have poles along $M''$. The restriction $\ol{\phi}|_{M''}$ depends only on $\phi$.\el

\bd{def7}The restriction $\ol{\phi}|_{M''}$ is called the \emph{residue} of $\phi$ and denoted $\mathrm{res}(\phi)$.\ed

The reader can check that the map $\phi\mapsto \mathrm{res}(\phi)$ defines a homomorphism from $R(\A_i)$ to $R(\A'')$. 
If $\phi\in R(\A_{i-1})$ then one can take $\ol{\phi} = 0$ and $\theta = \phi$ then the residue $\mathrm{res}(\phi) = 0$.
On the other hand $\mathrm{res}(\phi\wedge \psi_i) = \phi|_{K_i}$. 

\begin{theorem}[Leray's residue theorem] \label{thm4}The cohomology class of $\mathrm{res}(\phi)$ in $M''$ depends only on the cohomology class of $\phi$ in $M'$. \et 

The cohomology class of $\mathrm{res}(\phi)$ is called the \textit{residue class} and denoted $\mathrm{Res}(\phi)$. The mapping $\phi\mapsto \mathrm{Res}(\phi)$ gives us a homomorphism from $H^k(M')$ to $H^{k-1}(M'')$. The alternate construction of this homomorphism is the following composition
\begin{equation}
\mathrm{Res}\colon H^k(M_i)\to H^{k+1}(M_{i-1}, M_i)\to H^{k+1}(U, U_0)\to H^{k-1}(M'')\label{sec3eq1}
\end{equation}
where the first map is the connecting homomorphism in the long exact sequence of the pair $(M_{i-1}, M_i)$, the second map is the excision isomorphism and the last map is the (inverse) Thom isomorphism. 
As an immediate consequence of Theorem \ref{thm4} we get homomorphisms from a graded piece $R_k$ to $H^k(M)$ by sending forms to their residue classes. 
We encode this in the following commuting diagram:
\begin{equation}
\begin{CD} 
0 @>>> R_k(\A')    @>{\iota}>>   R_k(\A_i)  @>{\mathrm{res}}>>   R_{k-1}(\A'') @>>>  0\\
@.     @VVV @VVV @VVV @.\\
0 @>>> H^k(M') @>{\iota^*_{M_i}}>>  H^k(M_i)  @>{\mathrm{Res}}>>   H^{k-1}(M'')    @>>> 0.
\end{CD}\label{sec3eq2}
\end{equation}

\bl{lem5}
The vertical maps in the Diagram \ref{sec3eq2} are isomorphisms. 
\el

\bpr The proof is given by straightforward induction on the dimension of the ambient torus and the number of hypersurfaces in the arrangement. \epr

The following theorem is now immediate. 

\bt{thm5}
If $\A$ is a DR-type toric arrangement in an $l$-torus then the graded algebras $R(\A)$ and $H^*(M(\A))$ are isomorphic. \et 

\bpr
Consider the graded homomorphism from $R(\A)\to H^*(M(\A))$ which sends $\xi_i\mapsto [\xi_i]$ for $1\leq i\leq l$ and $\psi_j\mapsto [\psi_j]$ for $1\leq j\leq n$.
Since both the algebras are generated in degree $1$ using Lemma \ref{lem5} we see that this map induces an isomorphism of graded algebras.
\epr

Hence for DR-type toric arrangements the relations in $H^*(M(\A))$ are satisfied at the level of forms. We illustrate this in an example.

\be{ex7} Consider the arrangement from Examples \ref{ex4} and \ref{ex5}. We let $\xi_i$ denote the form $d\log z_i$ for $i = 1, 2$ and $\psi_i$ denote the form corresponding to the $i$th character in $\A$ for $1\leq i\leq 4$. The relations in $H^*(M(\A))$ are given by setting the following $2$-forms to zero 
\begin{multline*}
\{\xi_1\psi_1, \xi_2\psi_2, (\xi_1 + \xi_2)\psi_3, (\xi_1-\xi_2)\psi_4, \psi_1\psi_2 - \psi_1\psi_3 + \psi_2\psi_3 - \xi_2\psi_3,\\ \psi_1\psi_4 - \psi_1\psi_3 + \psi_2\psi_3 - \xi_2\psi_3- \psi_2\psi_4 -\xi_2\psi_1\}.
\end{multline*}
These relations can be found either by brute force or by using the isomorphism induced by the biholomorphism described in Example \ref{ex4}. The first $4$ relations arise due to dependency of the characters on the coordinate functions whereas the last two relations are due to the dependency in the characters themselves. It can be verified that $H^2(M(\A))$ is spanned by nine $2$-forms. Hence the Poincar\'e polynomial of the complement is $1 + 6t + 9t^2$.
\ee

\section{Concluding remarks}\label{sec5}
As pointed out by De Concini and Procesi in \cite[Remark 5.3]{copr05} the relations in $H^*(M(\A))$ for unimodular toric arrangements are quite complicated. We note here that the combinatorial description of the cohomology in case of DR-type arrangements will be treated in a forthcoming paper. In this last section we briefly outline the current work in progress.\par 

The reader will realize that in view of Remark \ref{rem2} the main theorem (Theorem \ref{thm1}) is true for any cohomology with complex coefficients. The logarithmic $1$-forms will have to be replaced by the pullbacks of a generator of $H^1(\C^*)$. Let $\A$ be a DR-type toric arrangement in $(\C^*)^l$. Consider the map $f\colon M(\A)\to (\C^*)^{l+n}$ whose first $l$ components are the coordinate functions $z_i$'s and the remaining $n$ components are the characters $\chi_i-c_i$ in $\A$. Theorem \ref{thm1} asserts that the induced map in cohomology $f^*\colon H^*((\C^*)^{l+n}) \to H^*(M(\A))$ is surjective. 
The ideal of relations $\ker f^*$ can be determined using deletion-restriction; this is work in progress.\par 

These relations will also help determine the conditions under which the (integer) cohomology is torsion free. For example, the relations described by De Concini and Procesi in case of unimodular arrangements \cite[Theorem 5.2]{copr05} have coefficients $\pm 1$, which are units in $\Z$, hence can be used to describe set of free generators for $H^*(M(\A); \Z)$. We would like to point the reader to the recent work by d'Antonio and Delucchi \cite{delu13}. With the help of discrete Morse theory the authors show that in case of the so-called complexified toric arrangements the cohomology of the complement is torsion free.\par

There is an interesting class of toric arrangements called \textit{toric Weyl arrangements}. Let $\Phi$ denote an irreducible, crystallographic root system. For every simple root $\alpha\in\Phi$ the function $\exp(\alpha)$ defines a character on an appropriate-dimensional torus. The collection of hypersurfaces corresponding to these characters is called the toric Weyl arrangement of type $\Phi$. 
In a forth-coming paper we show that these Weyl arrangements are of DR-type and compute the cohomology of the associated complement. 


\renewcommand{\bibname}{References}
\bibliographystyle{abbrv} 
\bibliography{bib_toric} 
\end{document}